\documentclass[a4paper,dvipdfmx]{amsart}
\let\oldtocsection=\tocsection
\let\oldtocsubsection=\tocsubsection
\let\oldtocsubsubsection=\tocsubsubsection
\renewcommand{\tocsection}[2]{\hspace{0em}\oldtocsection{#1}{#2}}
\renewcommand{\tocsubsection}[2]{\hspace{2em}\oldtocsubsection{#1}{#2}}
\renewcommand{\tocsubsubsection}[2]{\hspace{4em}\oldtocsubsubsection{#1}{#2}}
\usepackage{amssymb}
\usepackage{hyperref}
\usepackage{verbatim}
\usepackage{color}
\usepackage[dvipdfmx]{graphicx}
\theoremstyle{plain}
  \newtheorem{thm}{Theorem}[section]
  \newtheorem*{thm*}{Theorem}

  \newtheorem{conj}[thm]{Conjecture}

  \newtheorem*{obs*}{Observation}
\theoremstyle{definition}

\theoremstyle{remark}
  \newtheorem{rem}[thm]{Remark}

\newcommand{\Z}{\mathbb{Z}}
\newcommand{\C}{\mathbb{C}}

\newcommand{\Vol}{\operatorname{Vol}}
\newcommand{\CS}{\operatorname{CS}}

\newcommand{\Int}{\operatorname{Int}}

\newcommand{\Res}{\operatorname{Res}}

\newcommand{\Tr}{\operatorname{Tr}}
\newcommand{\SL}{\mathrm{SL}}
\renewcommand{\sl}{\mathfrak{sl}}

\newcommand{\AN}{\rm{AN}}
\newcommand{\NA}{\rm{NA}}
\newcommand{\NN}{\rm{NN}}
\numberwithin{equation}{section}
\allowdisplaybreaks

\begin{document}

\title{Kashaev invariants of twice-iterated torus knots}

\author{Hitoshi Murakami}
\address{
Graduate School of Information Sciences,
Tohoku University,
Aramaki-aza-Aoba 6-3-09, Aoba-ku,
Sendai 980-8579, Japan}
\email{starshea@tky3.3web.ne.jp}
\author{Anh T.~Tran}
\address{
Department of Mathematical Sciences, The University of Texas at Dallas, Richardson,
TX 75080, USA}
\email{att140830@utdallas.edu}
\date{\today}
\begin{abstract}
We calculate the asymptotic behavior of the Kashaev invariant of a twice-itarated torus knot and obtain topological interpretation of the formula in terms of the Chern--Simons invariant and the twisted Reidemeister torsion.
\end{abstract}
\keywords{volume conjecture; colored Jones polynomial; Kashaev invariant; iterated torus knot; Chern--Simons invariant; Reidemeister torsion}
\subjclass[2010]{Primary 57M27 57M25 57M50}
\thanks{H.M. was supported by JSPS KAKENHI Grant Numbers 26400079, 17K05239.
A.T. has been partially supported by a grant from the Simons Foundation (\#354595 to AT).}
\maketitle
\section{Introduction}
For a knot $K$ in the three-sphere $S^3$ and an integer $N\ge2$, let $\langle K\rangle_N\in\C$ be the Kashaev invariant \cite{Kashaev:MODPLA95}.
In \cite{Kashaev:LETMP97}, he conjectured that $|\langle K\rangle_N|$ grows exponentially with growth rate $\Vol(S^3\setminus{K})/(2\pi)$ for $N\to\infty$ when $K$ is hyperbolic, where $\Vol(S^3\setminus{K})$ is the hyperbolic volume of $S^3\setminus{K}$.
\par
J.~Murakami and the first author proved that the Kashaev invariant coincides with $J_N(K;\exp(2\pi i/N))$, where $J_N(K;q)$ is the colored Jones polynomial of a knot $K$ in the three-sphere $S^3$ associated with the $N$-dimensional irreducible representation of the Lie algebra $\sl_2(\C)$, normalized so that $J_N(U;q)=1$ for the unknot $U$ \cite{Murakami/Murakami:ACTAM12001}.
They also proposed the following conjecture:
\begin{conj}[Volume Conjecture]
For any knot $K$, we have
\begin{equation*}
  \frac{\log\left|J_N(K;\exp(2\pi i/N))\right|}{N}
  =
  \frac{\Vol(S^3\setminus{K})}{2\pi},
\end{equation*}
where $\Vol(S^3\setminus{K})$ is the simplicial volume of $S^3\setminus{K}$ normalized so that the simplicial volume of a hyperbolic knot complement equals its hyperbolic volume.
In particular, when $K$ is hyperbolic, Kashaev's conjecture holds.
\end{conj}
This conjecture was first proved for torus knots by Kashaev and O.~Tirkkonen \cite{Kashaev/Tirkkonen:ZAPNS2000}.
Note that since the complement of a torus knot is a Seifert fibered space, its simplicial volume is zero.
A proof of the conjecture for the figure-eight knot was given by T.~Ekholm (see for example \cite{Murakami/Yokota:2018} for the proof).
\par
J.~Andersen and S.~Hansen proved for the figure-eight knot $4_1$ the following asymptotic equivalence holds \cite{Andersen/Hansen:JKNOT2006}.
\begin{equation*}
  J_N(4_1;\exp(2\pi i/N))
  \underset{N\to\infty}{\sim}
  2\pi^{3/2}\left(\frac{N}{2\pi i}\right)^{3/2}\tau(4_1)
  \exp\left(\frac{N}{2\pi i}S(4_1)\right),
\end{equation*}
where $\tau(4_1):=\left(\frac{2}{\sqrt{-3}}\right)^{1/2}$ and $S(4_1):=\sqrt{-1}\Vol(S^3\setminus{4_1})$.
It is known that $\tau(4_1)^{-2}$ is the homological Reidemeister tosion twisted by the holonomy representation of $\pi_1(S^3\setminus{4_1})$ associated with the meridian, and $S(4_1)$ is the $\SL(2;\C)$ Reidemeister torsion of the holonomy representation.
It is also conjectured that for any hyperbolic knot $K$, we have
\begin{equation*}
  J_N(K;\exp(2\pi i/N))
  \underset{N\to\infty}{\sim}
  2\pi^{3/2}\left(\frac{N}{2\pi i}\right)^{3/2}\tau(K)
  \exp\left(\frac{N}{2\pi i}S(K)\right),
\end{equation*}
where $\tau(K)^{-2}$ and $S(K)$ are defined as above.
See \cite{Gukov/Murakami:FIC2008} and \cite{Dimofte/Gukov/Lenells/Zagier:CNTP2010}.
See also \cite{Ohtsuki:QT2016,Ohtsuki/Yokota:MATPC2018}.
\par
Let $T(c,d)$ be the torus knot of type $(c,d)$ for coprime integers $c$ and $d$.
J.~Dubois and Kashaev \cite{Dubois/Kashaev:MATHA2007} obtained the following formula:
\begin{equation}\label{eq:Dubois_Kashaev}
\begin{split}
  &J_N\bigl(T(c,d);e^{2\pi i/N}\bigr)
  \\
  \underset{N\to\infty}{\sim}&
  \frac{\pi^{3/2}}{2cd}
  \left(\frac{N}{2\pi i}\right)^{3/2}
  \sum_{k=1}^{cd-1}
  (-1)^{k+1}k^2
  \tau(k)
  \exp\left(S(k)\frac{N}{2\pi i}\right)
  +O(1),
\end{split}
\end{equation}
where
\begin{align*}
  S(k)&:=\frac{(k-cd)^2\pi^2}{cd},
  \\
  \tau(k)&:=\frac{4\sin(k\pi/c)\sin(k\pi/d)}{\sqrt{cd}}.
\end{align*}
See also \cite{Hikami/Murakami:Bonn} for the formulation above.
They also show that $S(k)$ is the Chern--Simons invariant and $\tau(k)^{-2}$ is the homological twisted Reidemeister torsion both associated with suitable irreducible representation from $\pi_1(S^3\setminus{T(c,d)})$ to $\SL(2,\C)$.
See also \cite{Kashaev/Tirkkonen:ZAPNS2000}.
\begin{rem}
If we define
\begin{equation*}
  \tilde{S}(k)
  :=
  \frac{k^2\pi^2}{cd},
\end{equation*}
the right hand side of \eqref{eq:Dubois_Kashaev} becomes
\begin{equation*}
  \frac{1}{2\sqrt{\pi}}
  \left(\frac{N}{2\pi i}\right)^{3/2}
  \sum_{k=1}^{cd-1}
  (-1)^{k(N+1)}i^{-cdN}
  \tau(k)\tilde{S}(k)
  \exp\left(\tilde{S}(k)\frac{N}{2\pi i}\right)
  +O(1).
\end{equation*}
Note that the Chern--Simons invariant is defined modulo $\pi^2\Z$ and that $S(k)\equiv\tilde{S}(k)\pmod{\pi^2\Z}$.
\end{rem}
\par
In \cite{Murakami:ACTMV2014} and \cite{Murakami:TOPOA2019}, the first author obtained a similar asymptotic formula for $J_N\bigl(T(2,2a+1)^{(2,2b+1)};\exp(\xi/N)\bigr)$ with $\xi\ne2\pi i$, where $T(2,2a+1)^{(2,2b+1)}$ is the $(2,2b+1)$-cable of $T(2,2a+1)$.
The purpose of this paper is to give an asymptotic formula for $J_N\bigl(T(2,2a+1)^{(2,2b+1)};\exp(2\pi i/N)\bigr).$
\begin{thm}
If $2b+1>4(2a+1)>0$, then we have
\begin{equation}\label{eq:J_N}
\begin{split}
  &J_N\bigl(T(2,2a+1)^{(2,2b+1)}, e^{2\pi i/N}\bigr)
  \\
  \underset{N\to\infty}{\sim}&
  \frac{1}{2\sqrt{\pi}}\left(\frac{N}{2\pi i}\right)^{3/2}
  \sum_{l=0}^{2b} \tau_1(l) S_1(l) \, e^{\frac{N}{2\pi i} S_1(l)}
  \\
  &+
  (-1)^{N}
  \frac{1}{2\sqrt{\pi}}\left(\frac{N}{2\pi i}\right)^{3/2}
  \sum_{m=0}^{4a+1} \tau_2(m) S_2(m) \, e^{\frac{N}{2\pi i} S_2(m)}
  \\
  &-
  \frac{1}{2}\left(\frac{N}{2\pi i}\right)^{2}
  \sum_{(j,k) \in \mathcal{B}} \tau_3(j,k) S_3(j,k) \, e^{\frac{N}{2\pi i} S_3(j,k)}
  + O(N^{1/2}),
\end{split}
\end{equation}
where 
\begin{align*}
  \tau_1(l)
  &=(-1)^l \sqrt{\frac{2}{2b+1}} \,
    \frac{\sin\left(\frac{2(2l+1)\pi}{2b+1}\right)}{\cos\left(\frac{(2a+1)(2l+1)\pi}{2b+1}\right)},
  \\
  S_1(l)
  &=\frac{(2l+1)^2\pi^2}{2(2b+1)},
  \\
  \tau_2(m)
  &=(-1)^m \sqrt{\frac{2}{2a+1}} \,
    \sin\left(\frac{(2m+1)\pi}{2a+1}\right),
  \\
  S_2(m)
  &=\frac{(2m+1)^2\pi^2}{2(2a+1)},
  \\
  \tau_3(j,k)
  &=(-1)^{j+k}\frac{4}{\sqrt{(2a+1)(2b+1-4(2a+1))}}
    \sin\left(\frac{(2k+1)\pi}{2a+1}\right),
  \\
  S_3(j,k)
  &=\left( \frac{(2k+1)^2}{2(2a+1)}+\frac{(2j+1)^2}{2\bigl(2b+1-4(2a+1)\bigr)} \right)\pi^2,
\end{align*}
and $\mathcal{B}$ is the set of all pairs of integers $(j,k)$ such that $0 \le k \le 4a+1$, $0 \le j \le 2b-4(2a+1)$, and $(2b+1-4(2a+1))(2k+1)<2(2a+1)(2j+1)$.
\end{thm}
We can also prove that $\tau_1(l)^{-2}$, $\tau_2(m)^{-2}$, and $\tau_3(l,m)^{-2}$ are the homological twisted Reidemeister torsions of certain representations of the fundamental group to $\rm{SL}(2;\C)$, and that $S_1(l)$, $S_2(m)$, and $S_3(l,m)$ are the Chern--Simons invariants of these representations.
\section{Proof of the asymptotic formula}
We will follow \cite[Section~5.2]{Liu:2008}.
Let $e_n$ be the $n$-th Jones--Wenzl idempotent in the Kauffman bracket skein algebra of an annulus defined by $z_1e_n=e_{n-1}+e_{n+1}$ with $e_0=1, e_1=z_1$, where $z_1$ is the circle around the annulus.
Let $\langle e_n\rangle_{\mathcal{K}^{\sigma}}$ be the Kauffman bracket of the element obtained from a framed knot $\mathcal{K}^{\sigma}$ by replacing a diagram of $\mathcal{K}^{\sigma}$ with $e_n$, where $\sigma$ is the framing.
If $\mathcal{T}^{2(2b+1)}$ is the $(2,2b+1)$-cable of the torus knot $T(2,2a+1)$ with framing $2(2b+1)$, then in \cite[Proposition~4]{Liu:2008} Q.~Liu proved
\begin{equation*}
\begin{split}
  &
  (-1)^{n}
  \exp\left(\frac{(2b+1-3(2a+1))\pi i}{4(n+1)}\right)
  (-n-1)
  \left(\frac{\langle e_n\rangle_{\mathcal{T}^{2(2b+1)}}}{[n+1]}\right)\Biggm|_{A=\exp(\pi i/(2(n+1)))}
  \\
  =&
  \frac{2(n+1)^3i}{\pi^4\sqrt{(2a+1)(2b+1-4(2a+1))}}
  \exp\left(\frac{-\pi i}{(2a+1)(n+1)}\right)
  \\&
  \times
  \int_{C_{\pi/4}}\int_{C_{\pi/4}}dz_1dz_2
  \psi_1(z_1)\psi_2(z_2) \delta(z_1,z_2) \, e^{(n+1)\theta (z_1,z_2)},
\end{split}
\end{equation*}
where
\begin{align*}
  [k]
  &:=\frac{A^{2k}-A^{-2k}}{A^2-A^{-2}},
  \\
  \theta (z_1,z_2)
  &:=
  -\frac{z_2^2}{(2a+1) \pi i}-\frac{4(z_1-z_2)^2}{(2b+1-4(2a+1)) \pi i}+4z_1,
  \\
  \delta(z_1, z_2)
  &:=
  \frac{z_2^2}{2a+1}+\frac{4(z_1-z_2)^2}{2b+1-4(2a+1)},
  \\
  \psi_1(z_1)
  &:=
  \frac{1}{2\cosh(2z_1)},
  \\
  \psi_2(z_2)
  &:=
  \frac{\sinh(\frac{2z_2}{2a+1})}{2\cosh(z_2)}.
\end{align*}
\par
Since the colored Jones polynomial $J_N(K;q)$ is normalized so that $J_N(U;q)=1$ with $U$ the unknot, if $K$ is the knot obtained from a framed knot $\mathcal{K}^{\sigma}$ by forgetting the framing, we have
\begin{equation*}
  J_N(K;q)
  =
  \frac{\langle e_{N-1}\rangle_{\mathcal{K}^{0}}}{\langle e_{N-1}\rangle_{\mathcal{U}^{0}}}
  =
  \bigl((-1)^{N-1}A^{N^2-1}\bigr)^{-\sigma}
  \frac{\langle e_{N-1}\rangle_{\mathcal{K}^{\sigma}}}{(-1)^{N+1}[N]},
\end{equation*}
where $\mathcal{K}^{0}$ is the $0$-framed knot obtained from $\mathcal{K}^{\sigma}$ by changing the framing and $\mathcal{U}^{0}$ is the framed unknot with framing $0$.
Therefore we have
\begin{equation*}
\begin{split}
  &J_N\bigl(T(2,2a+1)^{(2,2b+1)};e^{2\pi i/N}\bigr)
  \\
  =&
  (-1)^{N-1}
  \frac{1}{-N}
  \exp\left(-\frac{(2b+1-3(2a+1))\pi i}{4N}\right)
  \exp\left(\frac{-\pi i}{(2a+1)N}\right)
  \\
  &
  \left((-1)^{N-1}\exp\left(\frac{(N^2-1)\pi i}{2N}\right)\right)^{-2(2b+1)}
  (-1)^{N+1}
  \\
  &\times
  \frac{2N^3i}{\pi^4\sqrt{(2a+1)(2b+1-4(2a+1))}}
  \\&
  \times
  \int_{C_{\pi/4}}\int_{C_{\pi/4}}dz_1dz_2
  \psi_1(z_1)\psi_2(z_2) \delta(z_1,z_2) \, e^{N\theta (z_1,z_2)}
  \\
  =&
  (-1)^{N+1}
  e^{-\frac{\pi i}{4N}(2b+1-3(2a+1)+\frac{4}{2a+1})}
  e^{\frac{(2b+1)\pi i}{N}}
  \\
  &\times
  \left(\frac{2N^2 i}{\sqrt{(2a+1)(2b+1-4(2a+1))} \, \pi^4 } \right)
  \\&
  \times
  \int_{C_{\pi/4}}\int_{C_{\pi/4}}dz_1dz_2\psi_1(z_1)\psi_2(z_2) \delta(z_1,z_2) \, e^{N\theta (z_1,z_2)},
\end{split}
\end{equation*}
where
\begin{align*}
  \theta (z_1,z_2)
  &:=
  -\frac{z_2^2}{(2a+1) \pi i}-\frac{4(z_1-z_2)^2}{(2b+1-4(2a+1)) \pi i}+4z_1,
  \\
  \delta(z_1, z_2)
  &:=
  \frac{z_2^2}{2a+1}+\frac{4(z_1-z_2)^2}{2b+1-4(2a+1)},
  \\
  \psi_1(z_1)
  &:=
  \frac{1}{2\cosh(2z_1)},
  \\
  \psi_2(z_2)
  &:=
  \frac{\sinh(\frac{2z_2}{2a+1})}{2\cosh(z_2)}.
\end{align*}
\begin{rem}
We need to multiply by $(-1)^{N-1}\left((-1)^{N-1}\exp(\frac{(N^2-1)\pi i}{2N})\right)^{-2(2b+1)}$.
The first $(-1)^{N-1}$ is because Liu normalized the colored Jones polynomial by dividing the Kauffman bracket by $\frac{q^{N/2}-q^{-N/2}}{q^{1/2}-q^{-1/2}}$ but we need to divide it by $(-1)^{N-1}\frac{q^{N/2}-q^{-N/2}}{q^{1/2}-q^{-1/2}}$, which is the Kauffman bracket for the unknot.
The second one is because Liu's formula is for a framed knot with framing $2(2b+1)$.
\end{rem}
\begin{rem}
\begin{enumerate}
\item
$\theta(z_1,z_2)$ has a unique critical point
\begin{equation*}
  (w_1,w_2)=  \left( \frac{(2b+1)\pi i}2, 2(2a+1)\pi i \right).
\end{equation*}
\item
The poles of $\psi_1(z_1)$ between $C_{\pi/4}$ and $C_{\pi/4}+w_1$ are $\xi_l=\frac{2l+1}{4}\pi i$ ($0 \le l \le 2b$), where $C_{\pi/4}+w_1$ is the line passing through $w_1$ that is parallel to $C_{\pi/4}$.
Moreover, we have
\begin{equation*}
  \Res(\psi_1, \xi_l) = (-1)^{l-1} \, \frac{i}{4}.
\end{equation*}
\item
The poles of $\psi_2(z_2)$ between $C_{\pi/4}$ and $C_{\pi/4}+w_2$ are $\eta_m=\frac{2m+1}{2}\pi i$ ($0 \le m \le 4a+1$), where $C_{\pi/4}+w_2$ is the line passing through $w_2$ that is parallel to $C_{\pi/4}$.
Moreover, we have
\begin{equation*}
  \Res(\psi_2, \eta_m) = \frac{(-1)^{m}}{2} \sin(\frac{(2m+1)\pi}{2a+1}).
\end{equation*}
\end{enumerate}
\end{rem}
\par
Put
\begin{align*}
  F_N(z_1, z_2)
  &:=
  \delta(z_1,z_2) \, e^{N\theta (z_1,z_2)},
  \\
  I_N
  &:=
  \int_{C_{\pi/4}}\int_{C_{\pi/4}}dz_1dz_2\psi_1(z_1)\psi_2(z_2)F_N(z_1, z_2),
  \\
  \mathcal{J}_N
  &:=
  (-1)^{N-1} \left(\frac{2N^2 i}{\sqrt{(2a+1)(2b+1-4(2a+1))} \, \pi^4 } \right) I_N.
\end{align*}
Then $e^{-\frac{\pi i}{4N} \left( 2b+1-3(2a+1) + \frac{4}{2a+1} \right)} \, J_N(K, e^{2\pi i/N}) = \mathcal{J}_N $.
\par
By shifting the paths of integrations from $C_{\pi/4}$ to $C_{\pi/4}+w_1$ for the integration with respect to $z_1$, and from $C_{\pi/4}$ to $C_{\pi/4}+w_2$ for the integration with respect to $z_2$, we have
\begin{equation*}
  I_N
  =
  I_{0,N} + I_{1,N} + I_{2,N} + I_{3,N}
\end{equation*}
with
\begin{align*}
  I_{0,N}
  &:=
  \int_{C_{\pi/4} + w_1}\int_{C_{\pi/4} + w_2}\,dz_1\,dz_2 \psi_1(z_1)\psi_2(z_2)F_N(z_1,z_2),
  \\
  I_{1,N}
  &:=
  \sum_{l=0}^{2b} (2\pi i)\Res(\psi_1,\xi_l)\int_{C_{\pi/4}+w_2}\,dz_2\psi_2(z_2)F_N(\xi_l,z_2),
  \\
  I_{2,N}
  &:=
  \sum_{m=0}^{4a+1} (2\pi i)\Res(\psi_2,\eta_m)\int_{C_{\pi/4}+w_1}\,dz_1\psi_1(z_1)F_N(z_1,\eta_m),
  \\
  I_{3,N}
  &:=
  \sum_{l=0}^{2b}\sum_{m=0}^{4a+1}(2\pi i)^2\Res(\psi_1,\xi_l)\Res(\psi_2,\eta_m)F_N(\xi_l,\eta_m).
\end{align*}
Put $\mathcal{J}_{k,N} := (-1)^{N-1} \left(\frac{2N^2 i}{\sqrt{(2a+1)(2b+1-4(2a+1))} \, \pi^4 } \right) I_{k,N}$ so that
\begin{equation*}
  \mathcal{J}_{N}
  =
  \mathcal{J}_{0,N} + \mathcal{J}_{1,N} + \mathcal{J}_{2,N} + \mathcal{J}_{3,N}.
\end{equation*}
\subsection{\texorpdfstring{$I_{3,N}$}{I_{3,N}}}
In this subsection we calculate $I_{3,N}$.
\par
Since $\xi_l=\frac{2l+1}{4}\pi i$ and $\eta_m=\frac{2m+1}{2}\pi i$, we have 
\begin{equation*}
\begin{split}
  \delta(\xi_l,\eta_m)
  &=
  \frac{\eta_m^2}{2a+1}+\frac{4(\xi_l-\eta_m)^2}{2b+1-4(2a+1)}
  \\
  &= (\pi i)^2 \left( \frac{(2m+1)^2}{4(2a+1)}+\frac{(\frac{2l+1}{2}-(2m+1))^2}{2b+1-4(2a+1)} \right)
  =
  -\frac{1}{2}\tilde{S}_3(l,m),
\end{split}
\end{equation*}
where we put
\begin{equation*}
  \tilde{S}_3(l,m)
  :=
  2\pi^2 \left( \frac{(2m+1)^2}{4(2a+1)}+\frac{((2l+1)-2(2m+1))^2}{4(2b+1-4(2a+1))} \right).
\end{equation*}
Hence $\theta(\xi_l,\eta_m) = - (\pi i)^{-1} \delta(\xi_l,\eta_m) + 4 \xi_l = \frac{1}{2\pi i}\tilde{S}_3(l,m) + (2l+1) \pi i $ and
\begin{equation}\label{eq:S3}
  F_N(\xi_l,\eta_m)
  =
  \delta(\xi_l,\eta_m) \, e^{N\theta (\xi_l,\eta_m)}
  =
  \frac{(-1)^{N-1}}{2}\tilde{S}_3(l,m) \, e^{\frac{N}{2\pi i}\tilde{S}_3(l,m)}.
\end{equation}
Since $\Res(\psi_1,\xi_l)\Res(\psi_2,\eta_m) = (-1)^{l+m-1} \, \frac{i}{8} \sin(\frac{(2m+1)\pi}{2a+1})$, we obtain 
\begin{equation*}
  I_{3,N}
  =
  \sum_{l=0}^{2b}\sum_{m=0}^{4a+1}
  \pi^2 (-1)^{l+m+N-1} \, \frac{i}{4} \sin(\frac{(2m+1)\pi}{2a+1})\tilde{S}_3(l,m) \, e^{\frac{N}{2\pi i}\tilde{S}_3(l,m)}.
\end{equation*}
Hence 
\begin{equation*}
\begin{split}
  \mathcal{J}_{3,N}
  &=
  -\frac{N^2}{2\pi^2}
  \sum_{l=0}^{2b}\sum_{m=0}^{4a+1}
  \frac{(-1)^{l+m}}{\sqrt{(2a+1)(2b+1-4(2a+1))}}
  \\
  &\phantom{=-\frac{N^2}{2\pi^2}\sum_{l=0}^{2b}\sum_{m=0}^{4a+1}}
  \times\sin\left(\frac{(2m+1)\pi}{2a+1}\right)
  \tilde{S}_3(l,m)
  e^{\frac{N}{2\pi i}\tilde{S}_3(l,m)}
  \\
  &=
  - \frac{N^2}{8\pi^2} \sum_{l=0}^{2b} \sum_{m=0}^{4a+1} \tau_3(l,m)\tilde{S}_3(l,m) \, e^{\frac{N}{2\pi i}\tilde{S}_3(l,m)}.
\end{split}
\end{equation*}
\subsection{\texorpdfstring{$I_{2,N}$}{I_{2,N}}}
Next we calculate $I_{2,N}$.
\par
Consider the integral
\begin{equation*}
  \int_{C_{\pi/4}+w_1}\,dz_1\psi_1(z_1)F_N(z_1,\eta_m).
\end{equation*}
The polynomial $\theta (z_1,\eta_m) = -\frac{\eta_m^2}{(2a+1) \pi i}-\frac{4(z_1-\eta_m)^2}{(2b+1-4(2a+1)) \pi i}+4z_1$ has a unique critical point
\begin{equation*}
  \zeta_m:=\frac{2b+1-4(2a+1)}{2}\pi i + \eta_m= \frac{(2b+1) - 4(2a+1)+(2m+1)}{2}\pi i
\end{equation*}
and
\begin{equation*}
  \theta (z_1,\eta_m)=\theta (\zeta_m,\eta_m)-\frac{4}{(2b+1-4(2a+1))\pi i} (z_1-\zeta_m)^2.
\end{equation*}
Note that $\zeta_m$ is not a pole of $\psi_1(z_1)$ for any integer $m$.
\par
We have
\begin{equation*}
\begin{split}
  &\int_{C_{\pi/4}+w_1}\,dz_1\psi_1(z_1)F_N(z_1,\eta_m)
  \\
  =&
  \int_{C_{\pi/4}+\zeta_m}\,dz_1\psi_1(z_1)F_N(z_1,\eta_m)
  \\
  &-
  \sum_{l'=(2b+1) - 4(2a+1)+(2m+1) }^{2b} (2\pi i)\Res(\psi_1,\xi_{l'})F_N(\xi_{l'},\eta_m).
\end{split}
\end{equation*}
Note that the poles of $\psi(z_1)$ between $C_{\pi/4}+w_1$ and $C_{\pi/4}+\zeta_m$ are $\xi_{l'}=\frac{2l'+1}{r}\pi i$ ($(2b+1)-4(2a+1)+(2m+1)\le l'\le 2b$).
\par
Since $F_N(z_1,\eta_m)=e^{N \theta(\zeta_m,\eta_m)} \, \delta(z_1,\eta_m) \, e^{-\frac{4N}{(2b+1-4(2a+1))\pi i} (z_1-\zeta_m)^2}$ we have 
\begin{equation*}
\begin{split}
  &\int_{C_{\pi/4}+\zeta_m}\,dz_1\psi_1(z_1)F_N(z_1,\eta_m)
  \\
  =&
  e^{N\theta(\zeta_m,\eta_m)}\int_{C_{\pi/4}+\zeta_m}\,dz_1\psi_1(z_1)\delta(z_1,\eta_m)
  \, e^{-\frac{4N}{(2b+1-4(2a+1))\pi i} (z_1-\zeta_m)^2}
  \\
  =&
  e^{N\theta(\zeta_m,\eta_m)} \int_{C_{\pi/4}}\,dz_1   \psi_1(z_1+\zeta_m) \delta(z_1+\zeta_m, \eta_m)
  e^{-\frac {4N}{(2b+1-4(2a+1))\pi i}z_1^2}.
\end{split}
\end{equation*}
By the saddle point method (see, for example, \cite[Lemma~1]{Liu:2008}) we have
\begin{equation*}
\begin{split}
  &\int_{C_{\pi/4}}\,dz_1\psi_1(z_1+\zeta_m) \delta(z_1+\zeta_m, \eta_m) e^{-\frac {4N}{(2b+1-4(2a+1))\pi i}z_1^2}
  \\
  \underset{N\to\infty}{\sim}&
  \sqrt{\frac{(2b+1-4(2a+1)) \pi^2 i}{4N}} \, \psi_1(\zeta_m) \delta(\zeta_m, \eta_m)
  + O(N^{-3/2}).
\end{split}
\end{equation*}
Hence
\begin{equation*}
\begin{split}
  &\int_{C_{\pi/4}+\zeta_m}\,dz_1\psi_1(z_1)F_N(z_1,\eta_m)
  \\
  =&
  \sqrt{\frac{(2b+1-4(2a+1)) \pi^2 i}{4N}} \, \psi_1(\zeta_m)F_N(\zeta_m,\eta_m)
  + O(N^{-3/2}).
\end{split}
\end{equation*}
\par
Since $\zeta_m= \frac{(2b+1) - 4(2a+1)+(2m+1)}{2}\pi i$ and $\eta_m = \frac{2m+1}{2}\pi i$ we have
\begin{equation*}
\begin{split}
  \delta(\zeta_m,\eta_m)
  &=
  (\pi i)^2 \left( \frac{(2m+1)^2}{4(2a+1)}+2b+1-4(2a+1) \right)
  \\
  &=
  -\frac{1}{2} S_2(m)
  -\pi^2\bigl(2b+1-4(2a+1)\bigr).
\end{split}
\end{equation*}
Hence
\begin{equation*}
  \theta(\zeta_m,\eta_m)
  =
  - (\pi i)^{-1} \delta(\zeta_m,\eta_m) + 4\zeta_m
  =
  \frac{1}{2\pi i} S_2(m) + \bigl((2b+1) - 4(2a+1)+(2m+1)\bigr)\pi i
\end{equation*}
and 
\begin{equation*}
  F_N(\zeta_m,\eta_m)
  =
  -\frac{1}{2} S_2(m) \, e^{\frac{N}{2\pi i} S_2(m)}
  -\pi^2\bigl(2b+1-4(2a+1)\bigr)e^{\frac{N}{2\pi i}S_2(m)}.
\end{equation*}
\par
Since $\psi_1(\zeta_m) = \frac{1}{2\cosh (2\zeta_m)} = \frac{1}{2}$, we obtain
\begin{equation*}
\begin{split}
  &\int_{C_{\pi/4}+\zeta_m}\,dz_1\psi_1(z_1)F_N(z_1,\eta_m)
  \\
  =&
  -\frac{\pi}{8} \sqrt{\frac{(2b+1-4(2a+1)) i}{N}} S_2(m) \, e^{\frac{N}{2\pi i} S_2(m)}
  \\
  &-
  \frac{\pi^3}{4} \sqrt{\frac{(2b+1-4(2a+1)) i}{N}}\bigl((2b+1)-4(2a+1)\bigr)e^{\frac{N}{2\pi i} S_2(m)}
  \\
  &+ O(N^{-3/2}).
\end{split}
\end{equation*}
From \eqref{eq:S3} we have
\begin{equation*}
\begin{split}
  &I_{2,N}
  \\
  =&
  -
  \sum_{m=0}^{4a+1}\sum_{l'=(2b+1)-4(2a+1)+(2m+1)}^{2b}
  (2\pi i)^2 \Res(\psi_1,\xi_{l'})\Res(\psi_2,\eta_m)F_N(\xi_{l'},\eta_m)
  \\
  &
  +\sum_{m=0}^{4a+1}(2\pi i) \Res(\psi_2,\eta_m)\int_{C_{\pi/4}+\zeta_m}\,dz_1\psi_1(z_1)F_N(z_1,\eta_m)
  \\
  =&
  (-1)^{N+1}\frac{\pi^2i}{4}\sum_{m=0}^{4a+1}\sum_{l'=(2b+1)-4(2a+1)+(2m+1)}^{2b}
  (-1)^{l'+m-1}\sin\left(\frac{(2m+1)\pi}{2a+1}\right)
  \\
  &\phantom{(-1)^{N+1}\frac{\pi^2i}{4}\sum_{m=0}^{4a+1}\sum_{l'=(2b+1)-4(2a+1)+(2m+1)}^{2b}}
  \times
  \tilde{S}_3(l',m)e^{\frac{N}{2\pi i}\tilde{S}_3(l',m)}
  \\
  &-
  \frac{\pi^2i}{8}\sqrt{\frac{(2b+1-4(2a+1)) i}{N}}
  \sum_{m=0}^{4a+1}
  (-1)^{m}\sin\left(\frac{(2m+1)\pi}{2a+1}\right)S_2(m)e^{\frac{N}{2\pi i}S_2(m)}
  \\
  &-
  \frac{\pi^4i}{4} \sqrt{\frac{(2b+1-4(2a+1)) i}{N}}\bigl((2b+1)-4(2a+1)\bigr)
  \\
  &\quad
  \times
  \sum_{m=0}^{4a+1}(-1)^{m}\sin\left(\frac{(2m+1)\pi}{2a+1}\right)e^{\frac{N}{2\pi i} S_2(m)}
  \\
  &
  +O(N^{-3/2}).
\end{split}
\end{equation*}
Since $S_2(m+2a+1)=S_2(m)+4(a+m+1)\pi^2$, we have
\begin{equation*}
\begin{split}
  &\sum_{m=0}^{4a+1}(-1)^{m}\sin\left(\frac{(2m+1)\pi}{2a+1}\right)e^{\frac{N}{2\pi i} S_2(m)}
  \\
  =&
  \sum_{m=0}^{2a}(-1)^{m}\sin\left(\frac{(2m+1)\pi}{2a+1}\right)e^{\frac{N}{2\pi i} S_2(m)}
  \\
  &+
  \sum_{m=2a+1}^{4a+1}(-1)^{m}\sin\left(\frac{(2m+1)\pi}{2a+1}\right)e^{\frac{N}{2\pi i} S_2(m)}
  \\
  &\text{(put $m':=m-2a-1$ in the second term)}
  \\
  =&
  \sum_{m=0}^{2a}(-1)^{m}\sin\left(\frac{(2m+1)\pi}{2a+1}\right)e^{\frac{N}{2\pi i} S_2(m)}
  \\
  &+
  \sum_{m'=0}^{2a}(-1)^{m'+1}\sin\left(\frac{(2m'+1)\pi}{2a+1}\right)e^{\frac{N}{2\pi i} S_2(m')}
  \\
  =&0.
\end{split}
\end{equation*}
Therefore we have
\begin{equation*}
\begin{split}
  &I_{2,N}
  \\
  =&
  (-1)^{N+1}\frac{\pi^2i}{4}
  \sum_{m=0}^{4a+1}\sum_{l'=(2b+1)-4(2a+1)+(2m+1)}^{2b}
  (-1)^{l'+m-1}\sin\left(\frac{(2m+1)\pi}{2a+1}\right)
  \\
  &\phantom{(-1)^{N+1}\frac{\pi^2i}{4}\sum_{m=0}^{4a+1}\sum_{l'=(2b+1)-4(2a+1)+(2m+1)}^{2b}}
  \times
  \tilde{S}_3(l',m)e^{\frac{N}{2\pi i}\tilde{S}_3(l',m)}
  \\
  &-
  \frac{\pi^2i}{8}\sqrt{\frac{(2b+1-4(2a+1)) i}{N}}
  \sum_{m=0}^{4a+1}
  (-1)^{m}\sin\left(\frac{(2m+1)\pi}{2a+1}\right)S_2(m)e^{\frac{N}{2\pi i}S_2(m)}
  \\
  &
  +O(N^{-3/2}).
\end{split}
\end{equation*}
and so
\begin{equation*}
\begin{split}
  \mathcal{J}_{2,N}
  =&
  (-1)^{N-1} \left(\frac{2N^2i}{\sqrt{(2a+1)(2b+1-4(2a+1))} \, \pi^4} \right) I_{2,N}
  \\
  =&\frac{N^2}{8\pi^2}  \sum_{m=0}^{4a+1} \sum_{l'=(2b+1) - 4(2a+1)+(2m+1)}^{2b} \tau_3(l',m)\tilde{S}_3(l',m)
  \, e^{\frac{N}{2\pi i}\tilde{S}_3(l',m)}
  \\
  &+N^{3/2} (-1)^{N-1} \cdot \sqrt{\frac{i}{2}} \cdot \frac{ 1}{4\pi^2} \sum_{m=0}^{4a+1} \tau_2(m) S_2(m) \, e^{\frac{N}{2\pi i} S_2(m)} + O(N^{1/2}).
\end{split}
\end{equation*}
\subsection{\texorpdfstring{$I_{1,N}$}{I_{1,N}}}
Now we calculate $I_{1,N}$.
\par
Consider the integral
\begin{equation*}
  \int_{C_{\pi/4}+w_2}\,dz_2\psi_2(z_2)F_N(\xi_l,z_2).
\end{equation*}
The polynomial $\theta (\xi_l, z_2) = -\frac{z_2^2}{(2a+1) \pi i}-\frac{4(\xi_l-z_2)^2}{(2b+1-4(2a+1)) \pi i}+4\xi_l$ has a unique critical point
\begin{equation*}
  {\zeta'}_{l}:=\frac{4(2a+1)}{2b+1} \xi_l= \frac{(2a+1)(2l+1)}{2b+1}\pi i
\end{equation*}
and
\begin{equation*}
  \theta (\xi_l, z_2)=\theta (\xi_l,{\zeta'}_{l})-\frac{2b+1}{(2a+1)(2b+1-4(2a+1))\pi i} (z_2-{\zeta'}_{l})^2.
\end{equation*}
Note that ${\zeta'}_{l}$ is not a pole of $\psi_2(z_2)$ for any $l$.
\par
We have
\begin{equation*}
\begin{split}
  &\int_{C_{\pi/4}+w_2}\,dz_2\psi_2(z_2)F_N(\xi_l,z_2)
  \\
  =&
  \int_{C_{\pi/4}+{\zeta'}_{l}}\,dz_2\psi_2(z_2)F_N(\xi_l,z_2)
  -
  \sum_{m'=\lfloor \frac{(2a+1)(2l+1)}{2b+1} + \frac{1}{2} \rfloor }^{4a+1}
  (2\pi i)
  \Res(\psi_2,\eta_{m'})F_N(\xi_{l},\eta_{m'}).
\end{split}
\end{equation*}
Note that the pole of $\psi_2(z_2)$ between $C_{\pi/4}+w_2$ and $C_{\pi/4}+{\zeta'}_l$ are $\eta_{m'}$ ($\frac{(2a+1)(2l+1)}{2b+1}-\frac{1}{2}<m'<4a+\frac{3}{2}$, that is, $\lfloor \frac{(2a+1)(2l+1)}{2b+1} + \frac{1}{2} \rfloor\le m'\le4a+1$).
\par
By the same argument as in the previous case, we have
\begin{equation*}
\begin{split}
  &\int_{C_{\pi/4}+{\zeta'}_{l}}\,dz_2\psi_2(z_2)F_N(\xi_l,z_2)
  \\
  \underset{N\to\infty}{\sim}&
  \sqrt{\frac{(2a+1)(2b+1-4(2a+1)) \pi^2 i}{(2b+1)N}} \, \psi_2({\zeta'}_{l})F_N(\xi_l,{\zeta'}_{l})
  +
  O(N^{-3/2}).
\end{split}
\end{equation*}
Since $\xi_l= \frac{2l+1}{4}\pi i$ and ${\zeta'}_{l}= \frac{(2a+1)(2l+1)}{2b+1}\pi i$ we have
\begin{equation*}
  \delta(\xi_l,{\zeta'}_{l})
  = (\pi i)^2 \frac{(2l+1)^2}{4(2b+1)}
  = -\frac{1}{2} S_1(l).
\end{equation*}
Hence $\theta(\xi_l,{\zeta'}_{l}) = - (\pi i)^{-1} \delta(\xi_l,{\zeta'}_{l}) + 4 \xi_l = \frac{1}{2\pi i} S_1(l) + (2l+1)\pi i$ and
\begin{equation*}
  F_N(\xi_l,{\zeta'}_{l})
  =
  \frac{(-1)^{N-1}}{2} S_1(l) \, e^{\frac{N}{2\pi i} S_1(l)}.
\end{equation*}
Since
\begin{equation*}
  \psi_2({\zeta'}_{l})
  =
  \frac{\sinh(\frac{2{\zeta'}_{l}}{2a+1})}{2\cosh({\zeta'}_{l})}
  =
  i \frac{\sin(\frac{2(2l+1)}{2b+1}\pi)}{2\cos(\frac{(2a+1)(2l+1)}{2b+1}\pi)},
\end{equation*}
we obtain
\begin{equation*}
\begin{split}
  &\int_{C_{\pi/4}+{\zeta'}_{l}}\,dz_1\psi_1(z_1)F_N(z_1,\eta_m)
  \\
  =&
  (-1)^{N-1} \frac{\pi i}{4} \sqrt{\frac{(2a+1)(2b+1-4(2a+1))  i}{(2b+1)N}}
  \frac{\sin(\frac{2(2l+1)}{2b+1}\pi)}{\cos(\frac{(2a+1)(2l+1)}{2b+1}\pi)} S_1(l) \, e^{\frac{N}{2\pi i} S_1(l)}
  \\
  &+ O(N^{-3/2}).
\end{split}
\end{equation*}
\par
From \eqref{eq:S3} we have
\begin{equation*}
\begin{split}
  I_{1,N}
  \\
  =&
  -
  \sum_{l=0}^{2b}
  \sum_{m'=\lfloor \frac{(2a+1)(2l+1)}{2b+1} + \frac{1}{2} \rfloor }^{4a+1}
  (2\pi i)^2 \Res(\psi_1,\xi_l)\Res(\psi_2,\eta_{m'})F_N(\xi_l,\eta_{m'})
  \\
  &+
  \sum_{l=0}^{2b}
  (2\pi i) \Res(\psi_1,\xi_l)\int_{C_{\pi/4}+{\zeta'}_{l}}\,dz_2\psi_2(z_2)F_N(\xi_l,z_2)
  \\
  =&
  (-1)^{N+1}\frac{\pi^2i}{4}
  \sum_{l=0}^{2b}
  \sum_{m'=\lfloor\frac{(2a+1)(2l+1)}{2b+1}+\frac{1}{2}\rfloor}^{4a+1}
  (-1)^{l+m'-1}\sin(\frac{(2m'+1)\pi}{2a+1})
  \\
  &\phantom{(-1)^{N+1}\frac{\pi^2i}{4}\sum_{l=0}^{2b}
            \sum_{m'=\lfloor\frac{(2a+1)(2l+1)}{2b+1}+\frac{1}{2}\rfloor}^{4a+1}}
  \times
  \tilde{S}_3(l,m')e^{\frac{N}{2\pi i}\tilde{S}_3(l,m')}
  \\
  &+
  (-1)^{N}\frac{\pi^2i}{8}
  \sqrt{\frac{(2a+1)(2b+1-4(2a+1))i}{(2b+1)N}}
  \\
  &\quad\times
  \sum_{l=0}^{2b}
  (-1)^{l-1}
  \frac{\sin(\frac{2(2l+1)}{2b+1}\pi)}{\cos(\frac{(2a+1)(2l+1)}{2b+1}\pi)}S_1(l)e^{\frac{N}{2\pi i}S_1(l)}
  \\
  &+ O(N^{-3/2}).
\end{split}
\end{equation*}
Hence
\begin{equation*}
\begin{split}
  \mathcal{J}_{1,N}
  =&
  (-1)^{N-1} \left(\frac{2N^2i}{\sqrt{(2a+1)(2b+1-4(2a+1))} \, \pi^4} \right) I_{1,N}
  \\
  =&
  \frac{N^2}{8\pi^2}
  \sum_{l=0}^{2b}
  \sum_{m'=\lfloor \frac{(2a+1)(2l+1)}{2b+1}+\frac{1}{2}\rfloor}^{4a+1}
  \tau_3(l,m')\tilde{S}_3(l,m') \, e^{\frac{N}{2\pi i}\tilde{S}_3(l,m')}
  \\
  &-N^{3/2}  \cdot \sqrt{\frac{i}{2}} \cdot \frac{1}{4\pi^2}
  \sum_{l=0}^{2b} \tau_1(l) S_1(l)\,e^{\frac{N}{2\pi i} S_1(l)}
  + O(N^{1/2}).
\end{split}
\end{equation*}
\subsection{\texorpdfstring{$I_{0,N}$}{I_{0,N}}}
In this subsection we calculate $I_{0,N}$.
\par
We have
\begin{equation*}
\begin{split}
  &I_{0,N}
  \\
  =&
  \int_{C_{\pi/4} + w_1}\int_{C_{\pi/4} + w_2}\,dz_1\,dz_2 \psi_1(z_1)\psi_2(z_2) \delta(z_1,z_2)
  \, e^{N\theta (z_1,z_2)}
  \\
  =&
  \int_{C_{\pi/4}}\int_{C_{\pi/4}}\,dz_1\,dz_2 \psi_1(z_1+w_1)\psi_2(z_2+w_2) \delta(z_1+w_1,z_2+w_2)
  \, e^{N\theta (z_1+w_1,z_2+w_2)}.
\end{split}
\end{equation*}
Note that $\psi_1(z_1+w_1) = - \psi_1(z_1)$ and $\psi_2(z_2+w_2) =  \psi_2(z_2)$.
Moreover, we have
\begin{equation*}
\begin{split}
  \delta(z_1+w_1,z_2+w_2)
  &=
  \delta(z_1,z_2) + \delta(w_1,w_2) + \frac{2z_2w_2}{2a+1}+\frac{8(z_1-z_2)(w_1-w_2)}{2b+1-4(2a+1)}
  \\
  &=
  \delta(z_1,z_2) -\pi^2 (2b+1) + (4\pi i)z_1
\end{split}
\end{equation*}
and
\begin{equation*}
\begin{split}
  \theta(z_1+w_1,z_2+w_2)
  &=
  - (\pi i)^{-1} \delta(z_1+w_1,z_2+w_2) + 4(z_1+w_1)
  \\
  &=
  - (\pi i)^{-1}  \delta(z_1,z_2) - (2b+1) \pi i+ 4w_1
  \\
  &=
  - (\pi i)^{-1}  \delta(z_1,z_2) + (2b+1) \pi i.
\end{split}
\end{equation*}
Hence
\begin{equation*}
\begin{split}
  &I_{0,N}
  \\
  =&
  (-1)^{N-1}
  \\
  &\times
  \int_{C_{\pi/4}}\int_{C_{\pi/4}}\,dz_1\,dz_2 \psi_1(z_1)\psi_2(z_2)
  \big( \delta(z_1,z_2) -\pi^2 (2b+1) + (4\pi i)z_1 \big) \, e^{-\frac{N}{\pi i}\delta(z_1,z_2) }.
\end{split}
\end{equation*}
Since $\psi_1(z_1)\psi_2(z_2)$ is an odd function in $(z_1,z_2)$ and $\delta(z_1,z_2)$ is even, we obtain
\begin{equation*}
  I_{0,N}
  =
  (-1)^{N-1} \int_{C_{\pi/4}}\int_{C_{\pi/4}}\,dz_1\,dz_2 \psi_1(z_1)\psi_2(z_2) (4\pi i)z_1
  \, e^{-\frac{N}{\pi i}\delta(z_1,z_2) }.
\end{equation*}
From the saddle point method (see, for example, \cite[Lemma~3]{Liu:2008}), we have
\begin{equation*}
  I_{0,N}
  \underset{N\to\infty}{\sim}
  O(N^{-2}).
\end{equation*}
Hence 
\begin{equation*}
  \mathcal{J}_{0,N}
  =
  (-1)^{N-1} \left(\frac{2N^2i}{\sqrt{(2a+1)(2b+1-4(2a+1))} \, \pi^4} \right) I_{0,N}
  \underset{N\to\infty}{\sim}
  O(1).
\end{equation*}
\par
Finally we have
\begin{equation*}
\begin{split}
  J_N(K, e^{2\pi i/N})
  =&
  e^{\frac{\pi i}{4N}\left(2b+1-3(2a+1)+\frac{4}{2a+1}\right)}
  (\mathcal{J}_{0,N}+\mathcal{J}_{1,N}+\mathcal{J}_{2,N}+\mathcal{J}_{3,N})
  \\
  \underset{N\to\infty}{\sim}&
  -N^{3/2}  \cdot \sqrt{\frac{i}{2}} \cdot \frac{1}{4\pi^2}
  \sum_{l=0}^{2b} \tau_1(l) S_1(l)\,e^{\frac{N}{2\pi i} S_1(l)}
  \\
  &+
  N^{3/2} (-1)^{N-1} \cdot \sqrt{\frac{i}{2}} \cdot \frac{ 1}{4\pi^2}
  \sum_{m=0}^{4a+1} \tau_2(m) S_2(m) \, e^{\frac{N}{2\pi i} S_2(m)}
  \\
  &-
  \frac{N^2}{8\pi^2} \sum_{l=0}^{2b} \sum_{m=0}^{4a+1} \tau_3(l,m)\tilde{S}_3(l,m)
  \, e^{\frac{N}{2\pi i}\tilde{S}_3(l,m)}
  \\
  &+
  \frac{N^2}{8\pi^2}\sum_{m=0}^{4a+1}\sum_{l=(2b+1)-4(2a+1)+(2m+1)}^{2b} \tau_3(l,m)\tilde{S}_3(l,m)
  \, e^{\frac{N}{2\pi i}\tilde{S}_3(l,m)}
  \\
  &+
  \frac{N^2}{8\pi^2}
  \sum_{l=0}^{2b}
  \sum_{m=\lfloor \frac{(2a+1)(2l+1)}{2b+1}+\frac{1}{2}\rfloor}^{4a+1}
  \tau_3(l,m)\tilde{S}_3(l,m) \, e^{\frac{N}{2\pi i}\tilde{S}_3(l,m)}
  \\
  &+O(N^{1/2}).
\end{split}
\end{equation*}
The sum of the three double summations becomes
\begin{equation*}
  \sum_{(l,m)\in\mathcal{A}}\tau_3(l,m)\tilde{S}_3(l,m) \, e^{\frac{N}{2\pi i}\tilde{S}_3(l,m)},
\end{equation*}
where
\begin{equation*}
\begin{split}
  \mathcal{A}
  :=
  \{
  (l,m)\mid&
  0\le m \le 4a+1,
  0\le l\le 2b,
  l\le 2b+1-4(2a+1)+2m+1,
  \\
  &(2b+1)(2m+1)<2(2a+1)(2l+1)
  \}.
\end{split}
\end{equation*}
Putting $j:=l-(2m+1)$ and replacing $m$ with $k$, the summation above becomes
\begin{equation*}
\begin{split}
  &\sum_{(j,k)\in\mathcal{B}}
  \tau_3(j+2k+1,k)\tilde{S}_3(j+2k+1,j)e^{\frac{N}{2\pi i}\tilde{S}_3(j+2k+1,k)}
  \\
  =&
  -\sum_{(j,k)\in\mathcal{B}}
  \tau_3(j,k)S_3(j,k)e^{\frac{N}{2\pi i}S_3(j,k)},
\end{split}
\end{equation*}
where
\begin{equation*}
\begin{split}
  \mathcal{B}
  :=
  \{
  (j,k)\mid&
  0\le k\le 4a+1,
  0\le j\le 2b-4(2a+1),
  \\
  &(2b+1-4(2a+1))(2k+1)<2(2a+1)(2j+1)
  \}
\end{split}
\end{equation*}
and
\begin{equation*}
  S_3(j,k)
  :=
  2\pi^2\left(\frac{(2k+1)^2}{4(2a+1)}+\frac{(2j+1)^2}{4(2b+1-4(2a+1))}\right).
\end{equation*}
The theorem follows.
\section{Topological interpresentation}
In this section we give a toplogical interpretation of the right hand side of \eqref{eq:J_N}.
\subsection{Fundamental group}
We calculate the fundamental group of the complement of the twice-iterated torus knot $T(2,2a+1)^{(2,2b+1)}$.
\par
\begin{figure}[ht]
  \includegraphics[scale=0.3]{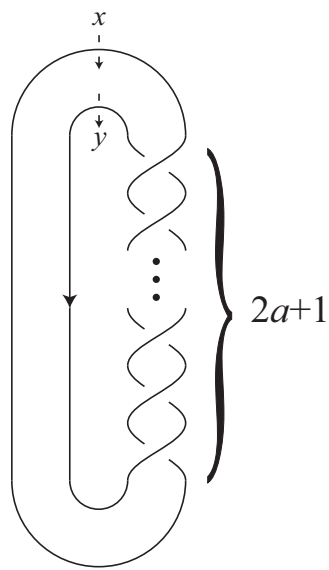}
  \caption{Torus knot $T(2,2a+1)$}
  \label{fig:torus_knot_pi1}
\end{figure}
\begin{figure}[ht]
  \includegraphics[scale=0.2]{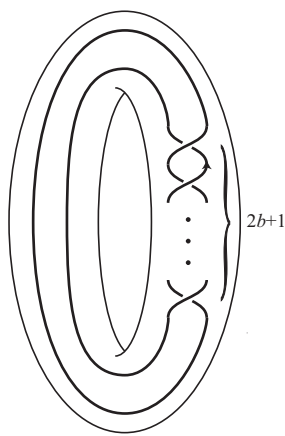}
  \caption{$(2,2b+1)$ torus knot in the solid torus}
  \label{fig:pattern}
\end{figure}

Put $E:=S^3\setminus\Int{N\left(T(2,2a+1)^{(2,2b+1)}\right)}$.
Then $E$ can be decomposed into $C:=S^3\setminus\Int N\bigl(T(2,2a+1)\bigr)$ (Figure~\ref{fig:torus_knot_pi1}) and $D$, where $N$ meas the regular neighborhood, $\Int$ means the interior, and $D$ is the complement of $(2,2b+1)$ torus knot in the solid torus (Figure~\ref{fig:pattern}).
Note that $D$ is homeomorphic to $S^3\setminus\Int N\bigl(T(2,4)\bigr)$ (see \eqref{eq:2_4}).
\begin{equation}\label{eq:2_4}
  \raisebox{-25mm}{\includegraphics[scale=0.2]{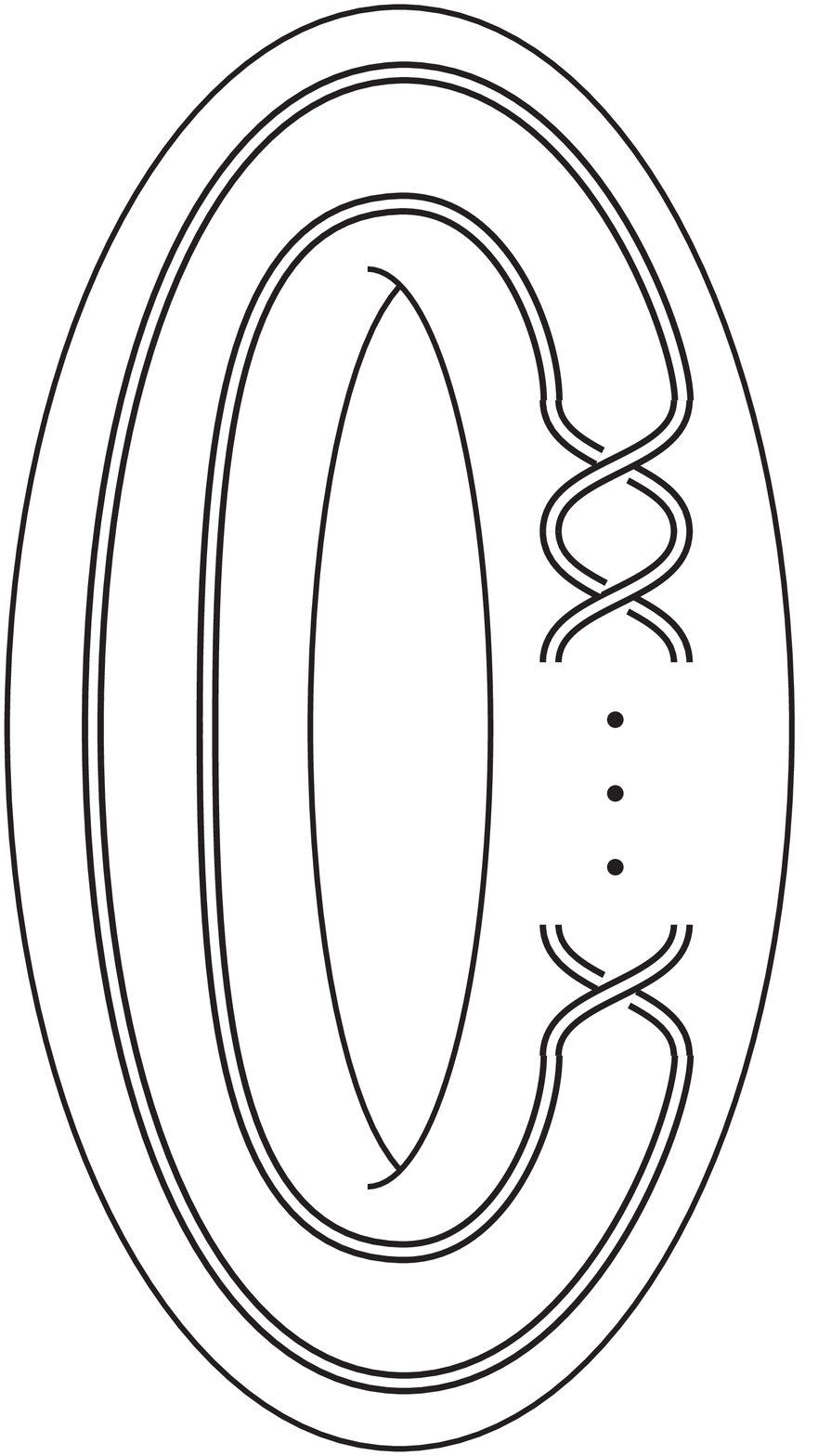}}
  \cong
  \raisebox{-25mm}{\includegraphics[scale=0.2]{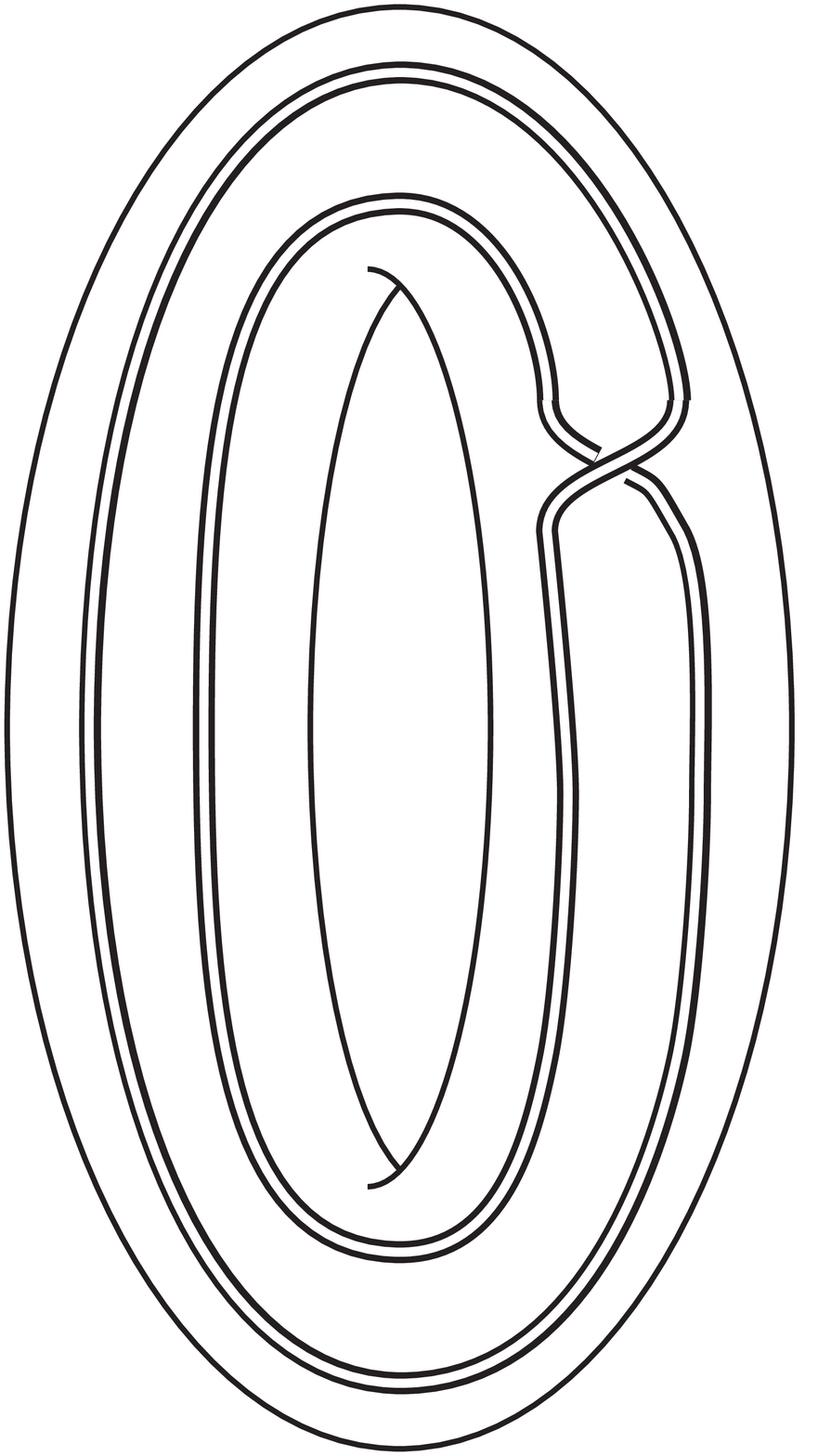}}
  \cong
  \raisebox{-25mm}{\includegraphics[scale=0.2]{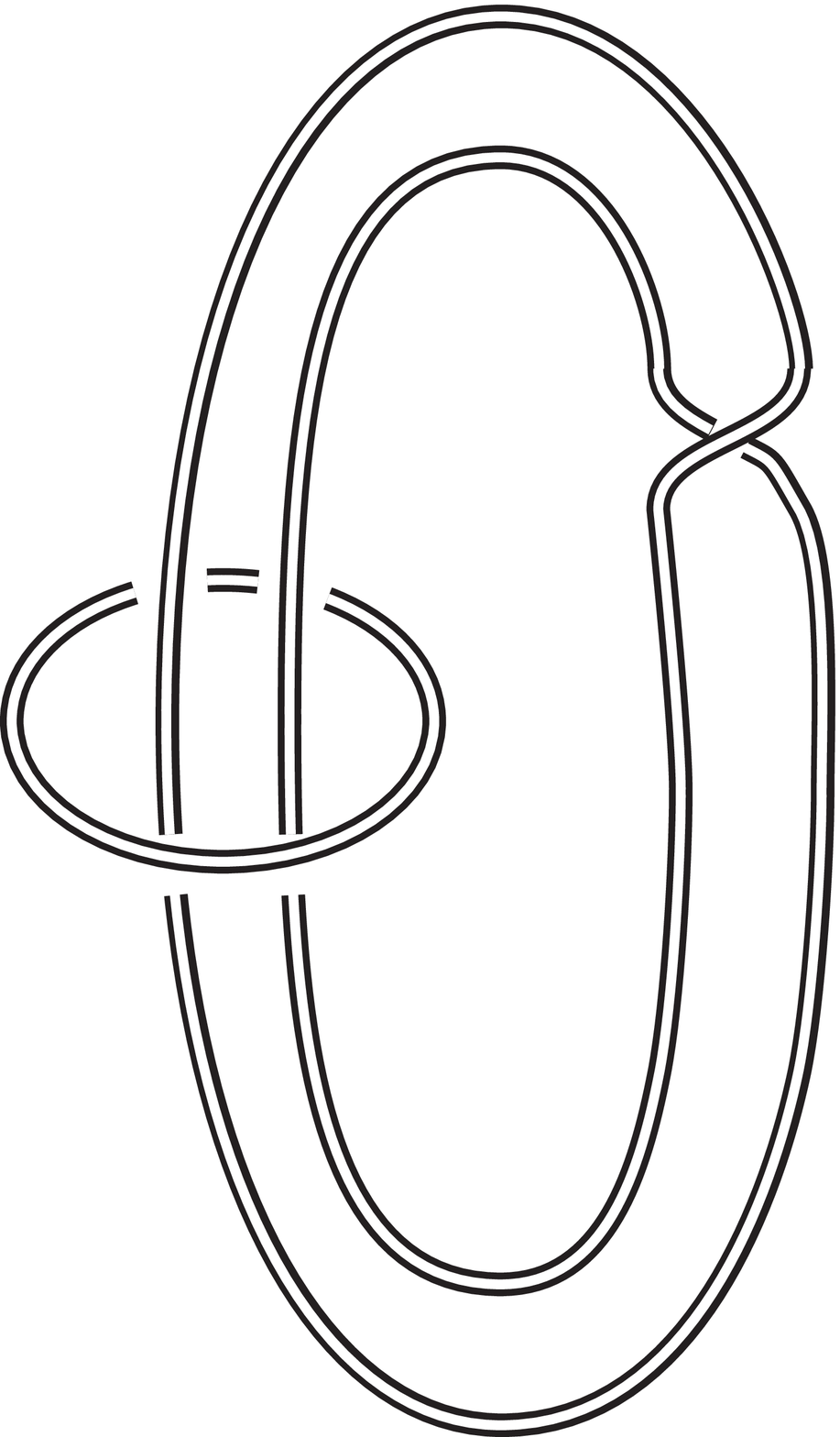}}
  \cong
  \raisebox{-18mm}{\includegraphics[scale=0.2]{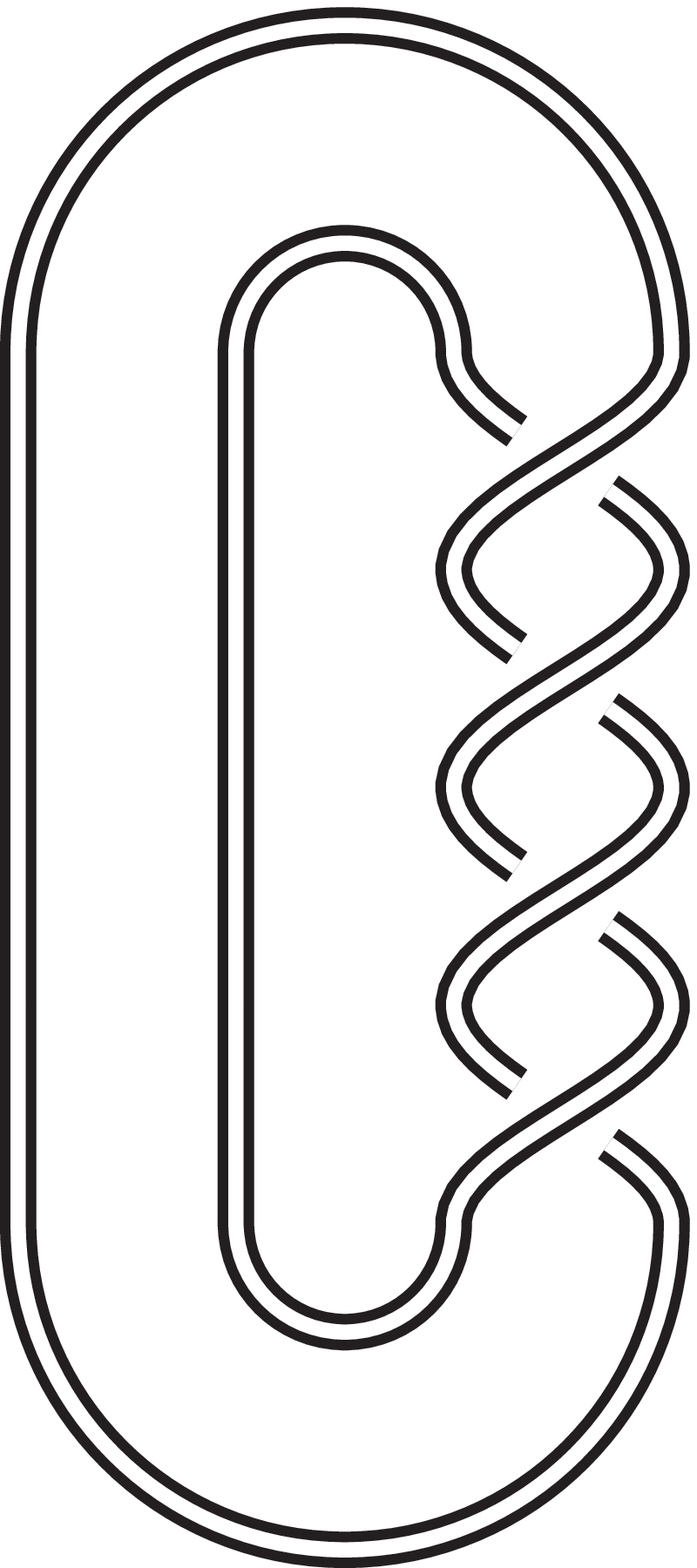}}
\end{equation}
\par
\begin{figure}[ht]
  \includegraphics[scale=0.3]{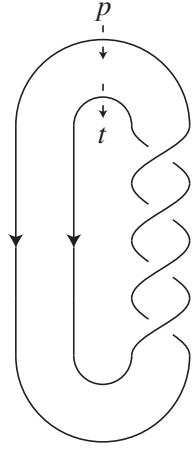}
  \caption{Torus link $(2,4)$}
  \label{fig:2_4}
\end{figure}
\begin{figure}[ht]
  \includegraphics[scale=0.2]{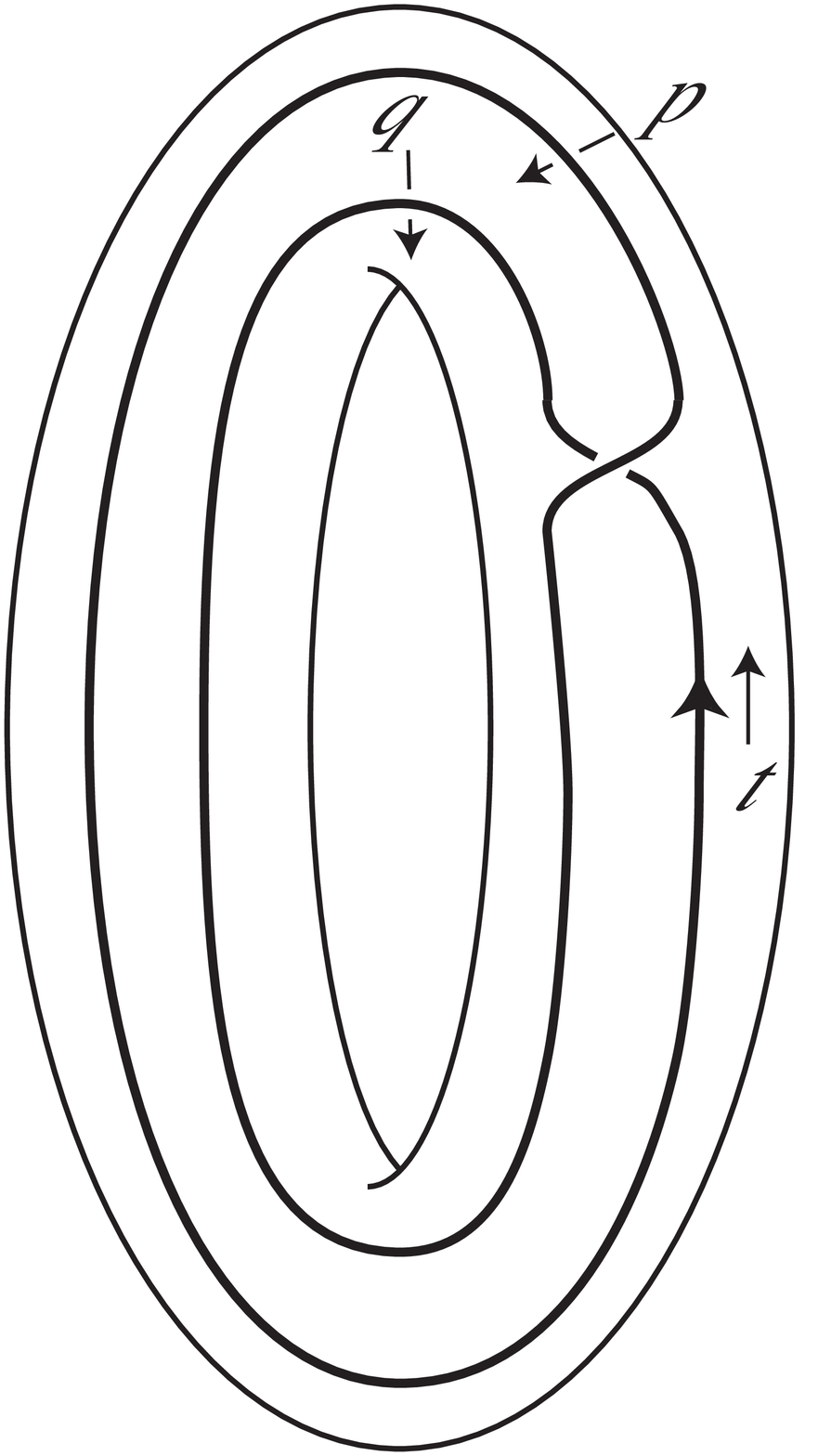}
  \caption{$q=tpt^{-1}$}
  \label{fig:pattern_untwisted_pi1}
\end{figure}

If we choose $x$ and $y$ as generators of $\pi_1(C)$ as indicated in Figure~\ref{fig:torus_knot_pi1}, we have
\begin{equation*}
  \pi_1(D)
  =
  \langle x,y\mid (xy)^ax=y(xy)^a\rangle,
\end{equation*}
where we choose the basepoint on the boundary of $N\bigl(T(2,2a+1)\bigr)$.
Let $p$ and $t$ be generators of $\pi_1\left(S^3\setminus\Int{N\bigl(T(2,4)\bigr)}\right)$ as in Figure~\ref{fig:2_4}.
Then we have
\begin{equation*}
  \pi_1\left(S^3\setminus\Int N\bigl(T(2,4)\bigr)\right)=\langle p,t\mid ptpt=tptp\rangle,
\end{equation*}
where the basepoint is at the bottom-right of the torus.
Let $q\in\pi_1(D)$ be the element as indicated in Figure~\ref{fig:pattern_untwisted_pi1}.
Then we see that $q=tpt^{-1}$
From van Kampen's theorem we have
\begin{align*}
  \pi_1(E)
  =
  \langle
    x,y,p,t\mid (xy)^ax=y(xy)^a,y(xy)^{2a}x^{-4a-1}=tx^{-b},x=ptpt^{-1}
  \rangle.
\end{align*}
Moreover we can see that the meridian $\mu$ and the preferred longitude $\lambda$ of $T(2,2a+1)^{(2,2b+1)}$ are given as follows:
\begin{align*}
  \mu
  &=p,
  \\
  \lambda
  &=
  y(xy)^{2a}x^{b-4a-1}p(tpt^{-1})^{-b}y(xy)^{2a}x^{b-4a-1}p^{-3b-1}.
\end{align*}
\subsection{Representation}\label{subsec:representation}
In this subsection, we construct non-Abelian representations from $\pi_1(E)$ to $\SL(2;\C)$.
Note that we do not know whether we exhaust all such representations.
\par
Put $\omega_1:=\exp\left(\frac{\pi i}{2b+1}\right)$, $\omega_2:=\exp\left(\frac{\pi i}{2a+1}\right)$, and $\omega_3:=\exp\left(\frac{\pi i}{2b+1-4(2a+1)}\right)$.
Let $\rho^{\rm{AN}}_{u;j}$ be the representation $\pi_1(E)\to\SL(2;\C)$ defined by
\begin{align*}
  \rho^{\rm{AN}}_{u;l}(x)
  &:=
  \rho^{\rm{AN}}_{u;l}(y)
  =
  T_{1}^{-1}
  \begin{pmatrix}
    \omega_1^{2l+1}&e^{-u/2} \\
    0       &\omega_1^{-(2l+1)}
  \end{pmatrix}
  T_{1},
  \\
  \rho^{\rm{AN}}_{u;l}(p)
  &:=
  \begin{pmatrix}
    e^{u/2}&1 \\
    0      &e^{-u/2}
  \end{pmatrix},
  \\
  \rho^{\rm{AN}}_{u;l}(t)
  &:=
  T_{1}^{-1}
  \begin{pmatrix}
    \omega_1^{(2l+1)b}&\frac{\omega_1^{(2l+1)b}-\omega_1^{-(2l+1)b}}{\omega_1^{2l+1}-\omega_1^{-(2l+1)}}e^{-u/2} \\
    0                 &\omega_1^{-(2l+1)b}
  \end{pmatrix}
  T_{1},
\end{align*}
where
\begin{equation*}
  T_{1}
  :=
  \begin{pmatrix}
    1                                 & 0 \\
    \omega_1^{-(2l+1)}e^{u/2}-e^{-u/2}&1
  \end{pmatrix}
\end{equation*}
and $l=0,1,\dots,b-1$.
Note that
\begin{equation*}
  \rho^{\AN}_{u;l+2b+1}=\rho^{\AN}_{u;l}
\end{equation*}
and that $\rho^{\AN}_{u,2b-l}$ is conjugate to $\rho^{\AN}_{u,l}$ by
\begin{equation*}
  R_{1,l}
  :=
  \begin{pmatrix}
    1                                          &0 \\
    e^{-u/2}-\omega_1^{2l+1}e^{u/2}&1
  \end{pmatrix}.
\end{equation*}
The longitude $\lambda$ is sent to
\begin{equation*}
  \rho^{\AN}_{u;l}(\lambda)
  =
  \begin{pmatrix}
    -e^{-(2b+1)u}& \frac{e^{(2b+1)u}-e^{-(2b+1)u}}{e^{u/2}-e^{-u/2}}\\
    0            &-e^{(2b+1)u}
  \end{pmatrix}.
\end{equation*}
\par
Let $\rho^{\rm{NA}}_{u;m}$ be the representation with
\begin{align*}
  \rho^{\rm{NA}}_{u;m}(x)
  &:=
  \begin{pmatrix}
    e^{u}&e^{u/2}+e^{-u/2} \\
    0    &e^{-u}
  \end{pmatrix},
  \\
  \rho^{\rm{NA}}_{u;m}(y)
  &:=
  \begin{pmatrix}
    e^{u}                                                                     &0 \\
    \frac{\omega_2^{2m+1}+\omega_2^{-(2m+1)}-e^{2u}-e^{-2u}}{e^{u/2}+e^{-u/2}}&e^{-u}
  \end{pmatrix},
  \\
  \rho^{\rm{NA}}_{u;m}(p)
  &:=
  \begin{pmatrix}
    e^{u/2}&1 \\
    0      &e^{-u/2}
  \end{pmatrix},
  \\
  \rho^{\rm{NA}}_{u;m}(t)
  &:=
  \begin{pmatrix}
    -e^{(b-4a-2)u}&\frac{e^{-(b-4a-2)u}-e^{(b-4a-2)u}}{e^{u/2}-e^{-u/2}}\\
    0             &-e^{-(b-4a-2)u}
  \end{pmatrix}
\end{align*}
with $m=0,1,\dots,a-1$.
The longitude $\lambda$ is sent to
\begin{equation*}
  \rho^{\NA}_{u;m}(\lambda)
  =
  \begin{pmatrix}
    e^{-4(2a+1)u}&\frac{e^{-4(2a+1)u}-e^{4(2a+1)u}}{e^{u/2}-e^{-u/2}} \\
    0            &e^{4(2a+1)u}
  \end{pmatrix}
  =
  \rho^{\NA}_{u;m}(p)^{-8(2a+1)}.
\end{equation*}
Note the following symmetries:
\begin{align*}
  \rho^{\rm{NA}}_{u;k+2a+1}&=\rho^{\rm{NA}}_{u;k},
  \\
  \rho^{\rm{NA}}_{u;2a-k}&=\rho^{\rm{NA}}_{u;k}.
\end{align*}
\par
Let $\rho^{\rm{NN}}_{u;j,k}$ be the representation with
\begin{align*}
  \rho^{\rm{NN}}_{u;j,k}(x)
  &:=
  T_2^{-1}
  \begin{pmatrix}
    \omega_3^{2j+1}&1 \\
    0              &\omega_3^{-(2j+1)}
  \end{pmatrix}
  T_2,
  \\
  \rho^{\rm{NN}}_{u;j,k}(y)
  &:=
  T_2^{-1}
  \begin{pmatrix}
    \omega_3^{2j+1}                                                          &0 \\
    \omega_2^{2k+1}+\omega_2^{-(2k+1)}-\omega_3^{2(2j+1)}-\omega_3^{-2(2j+1)}&\omega_3^{-(2j+1)}
  \end{pmatrix}
  T_2,
  \\
  \rho^{\rm{NN}}_{u;j,k}(p)
  &:=
  \begin{pmatrix}
    e^{u/2}&1 \\
    0      &e^{-u/2}
  \end{pmatrix},
  \\
  \rho^{\rm{NN}}_{u;j,k}(t)
  &:=
  T_2^{-1}
  \begin{pmatrix}
    -\omega_3^{(b-4a-2)(2j+1)}
      &\frac{\omega_3^{-(b-4a-2)(2j+1)}-\omega_3^{(b-4a-2)(2j+1)}}{\omega_3^{2j+1}-\omega_3^{-(2j+1)}}
    \\
    0 &-\omega_3^{-(b-4a-2)(2j+1)}
  \end{pmatrix}
  T_2,
\end{align*}
where
\begin{equation*}
  T_2
  :=
  \begin{pmatrix}
    e^{u/4}                            &0\\
    \omega_3^{-(2j+1)}e^{u/4}-e^{-3u/4}&e^{-u/4}
  \end{pmatrix},
\end{equation*}
$=0,1,\dots,2b-4(2a+1)$, and $k=0,1,\dots,a-1$.
The longitude $\lambda$ is sent to
\begin{equation*}
  \rho^{\rm{NN}}_{u;j,k}(\lambda)
  =
  \begin{pmatrix}
    -e^{-(2b+1)u}&\frac{e^{(2b+1)u}-e^{-(2b+1)u}}{e^{u/2}-e^{-u/2}} \\
    0            &-e^{(2b+1)u}
  \end{pmatrix}.
\end{equation*}
\par
Note the following symmetries:
\begin{align*}
  \rho^{\NN}_{u;j,k+2a+1}&=\rho^{\NN}_{u;j,k},
  \\
  \rho^{\NN}_{u;j,2a-k}&=\rho^{\NN}_{u;j,k},
  \\
  \rho^{\NN}_{u;j+2b+1-4(2a+1),k}&=\rho^{\NN}_{u;j,k},
\end{align*}
\begin{rem}
Note that $\rho^{\NN}_{0;j,k}$ and $\rho^{\NN}_{0;2b-4(2a+1)-j,k}$ are not conjugate because $\Tr\rho^{\NN}_{0;j,k}(py)=-(\omega_3^{2j+1}-1)^2+\omega_2^{2k+1}+\omega_2^{-(2k+1)}$.
\end{rem}
\subsection{Chern--Simons invariant}
For a knot $K$ in $S^3$, let $M$ be the complement of the interior of the regular neighborhood of $K$.
Denote by $\mu$ and $\lambda$ the meridian and the preferred longitude of $M$, respectively.
By a conjugation we may assume that a representation $\rho\colon M\to\SL(2;\Z)$ sends $\mu$ and $\lambda$ to
\begin{align*}
  \rho(\mu)&=\begin{pmatrix}e^{u/2}&\ast\\0&e^{-u/2}\end{pmatrix},
  \\
  \rho(\lambda)&=\begin{pmatrix}e^{v/2}&\ast\\0&e^{-v/2}\end{pmatrix},
\end{align*}
respectively.
The $\SL(2;\C)$ Chern--Simons invariant is a map from $X(M)$ to $\C$ modulo $\pi^2\Z$, where $X(M)$ is the $\SL(2;\C)$ character variety of $M$.
Note that we need to fix log branches of the eigenvalues $e^{u/2}$ and $e^{v/2}$.
See \cite{Kirk/Klassen:COMMP93} for details.
\par
In \cite[\S~2.5]{Murakami:ACTMV2014} the first author proved that the Chern--Simons invariants of $\rho^{\AN}_{u;l}$, $\rho^{\NA}_{u;m}$, $\rho^{\NN}_{u;j,k}$ are given as follows.
\begin{thm}[\cite{Murakami:ACTMV2014}]
Let $\rho^{\AN}_{u;l}$, $\rho^{\NA}_{u;m}$, and $\rho^{\NN}_{u;j,k}$ be representations given in Subsection~\ref{subsec:representation}.
Then we have
\begin{align*}
  \CS(\rho^{\rm{AN}}_{u;j})
  &=
  \frac{(2l+1)^2\pi^2}{2(2b+1)}+\frac{1}{2} d_1 u\pi i,
  \\
  \CS(\rho^{\rm{NA}}_{u;k})
  &=
  \frac{(2m+1)^2\pi^2}{2(2a+1)}+ d_2 u \pi i,
  \\
  \CS(\rho^{\rm{NN}}_{u;j,k})
  &=
  \frac{(2k+1)^2\pi^2}{2(2a+1)}
  +
  \frac{(2j+1)^2\pi^2}{2(2b+1-4(2a+1))}
  +
  \frac{1}{2}d_3u\pi i,
\end{align*}
for an odd integer $d_1$, and integers $d_2$ and $d_3$.
Here we choose $-2(2b+1)u+2d_1\pi i$, $-8(2a+1)u+4d_2\pi i$, and $-2(2b+1)u+2d_3\pi i$ as lifts of the eigenvalues of $\rho^{\AN}_{u;l}(\lambda)$, $\rho^{\NA}_{u;m}(\lambda)$, and $\rho^{\NN}_{u;j,k}(\lambda)$, respectively.
\end{thm}
Therefore we conclude
\begin{align*}
  S_1(l)
  &=
  \CS(\rho^{\rm{AN}}_{0;l}),
  \\
  S_2(m)
  &=
  \CS(\rho^{\rm{NA}}_{0;m}),
  \\
  S_3(l,m)
  &=
  \CS(\rho^{\rm{NN}}_{0;j,k}).
\end{align*}
\subsection{Reidemeister torsion}
Let $M$ be the complement of the interior of the regular neighborhood of a knot, and $\tilde{M}$ be its universal covering space.
Then the chain group $C_j(\tilde{M};\Z)$ can be regarded as a $\Z[\pi_1(M)]$-module.
Given a representation $\rho\colon\pi_1(M)\to\SL(2;\C)$, we can also regard the Lie algebra $\sl(2;\C)$ as a $\Z[\pi_1(M)]$-module by the adjoint action.
So we can define the tensor product $C_j(\tilde{M};\Z)\otimes_{\Z[\pi_1(M)]}\sl(2;\C)$ and denote it by $C_j(M;\rho)$.
The Reidemeister torsion of the corresponding chain complex $\{C_j(M;\rho),\partial_j\}$ is denoted by $\mathbb{T}(M;\rho)$ and called the homological twisted Reidemeister torsion of $M$ associated with $\rho$.
Note that we need to specify bases of $H_j(M;\rho)$ unless $\{C_j(M;\rho),\partial_j\}$ is acyclic, where $H_j(M;\rho)$ is the $j$-th homology group of the chain complex $\{C_j(M;\rho),\partial_j\}$.
\par
In \cite{Murakami:TOPOA2019}, the first author calculated the homplogical twisted Reidemeister torsions of $S^3\setminus\Int{N\left(T(2,2a+1)^{(2,2b+1)}\right)}$ associated with $\rho^{\AN}_{u,l}$, $\rho^{\NA}_{u;m}$ and $\rho^{\NN}_{u;j,k}$.
\begin{thm}[\cite{Murakami:TOPOA2019}]
Put $M:=S^3\setminus\Int{N\left(T(2,2a+1)^{(2,2b+1)}\right)}$.
The homological twisted Reidemeister torsions of $M$ associated with the representations defined in Subsection~\ref{subsec:representation} are given as follows.
\begin{align*}
  \mathbb{T}\left(M;\rho^{\AN}_{u;l}\right)
  &=
  \frac{(2b+1)\cos^2\left(\frac{(2a+1)(2l+1)\pi}{2b+1}\right)}{2\sin^2\left(\frac{2(2l+1)\pi}{2b+1}\right)},
  \\
  \mathbb{T}\left(M;\rho^{\NA}_{u;m}\right)
  &=
  \frac{(2a+1)}{2\sin^2\left(\frac{(2m+1)\pi}{2a+1}\right)},
  \\
  \mathbb{T}\left(M;\rho^{\NN}_{u;j,k}\right)
  &=
  \frac{(2a+1)(2b+1-4(2a+1))}{16\sin^2\left(\frac{(2k+1)\pi}{2a+1}\right)}.
\end{align*}
Note that since we have
\begin{align*}
  \dim H_j(M;\rho^{\AN})
  &=
  \begin{cases}
    1&\text{if $j=1,2$}
    \\
    0&\text{otherwise},
  \end{cases}
  \\
  \dim H_j(M;\rho^{\NA})
  &=
  \begin{cases}
    1&\text{if $j=1,2$}
    \\
    0&\text{otherwise},
  \end{cases}
  \\
  \dim H_j(M;\rho^{\NN})
  &=
  \begin{cases}
    2&\text{if $j=1,2$}
    \\
    0&\text{otherwise},
  \end{cases}
\end{align*}
we need to specify bases of non-trivial homology groups.
See \cite{Murakami:TOPOA2019} for details.
\end{thm}
\bibliography{mrabbrev,hitoshi}
\bibliographystyle{amsplain}
\end{document}